\documentclass[a4paper]{article}
\date{Nov 5, 2004}
\usepackage[T1]{fontenc}
\usepackage[latin1]{inputenc}
\usepackage{amsmath,amssymb,amsthm,paralist}
\usepackage{mathrsfs}
\usepackage{textcomp,xspace}

\newtheorem{theorem}{Theorem}
\newtheorem{lemma}[theorem]{Lemma}

\newtheorem{corollary}[theorem]{Corollary}

\newtheorem{fact}[theorem]{Fact}

\theoremstyle{definition}
\newtheorem{definition}[theorem]{Definition}
\newtheorem{notation}[theorem]{Notation}

\theoremstyle{remark}
\newtheorem*{remark}{Remark}

\catcode`\"\active
\newcommand{"}[2]{{}_{#1}{#2}}

\newcommand{\osum}[1]{#1_0+#1_1+#1_2\ldots}
\newcommand{\sumo}[1]{\ldots+"2{#1}+"1{#1}+"0{#1}}

\newcommand{\sC}{\mathscr{C}}
\let\hungarian\H
\renewcommand{\H}{\ensuremath{\sC'}}
\newcommand{\ncontemb}{\npreceq}
\newcommand{\contemb}{\preceq}

\newcommand{\dom}{\operatorname{\mathrm{dom}}}
\newcommand{\rk}{\operatorname{\mathrm{rk}}}
\newcommand{\cbrk}{\operatorname{\mathrm{rk}_{\text{CB}}}}

\newcommand{\gdl}[1]{{\bf G}_{#1}}

\newcommand{\OO}{{\mathscr O}}
\newcommand{\CB}{{\rm CB}}
\newcommand{\LL}{{\mathscr L}}

\newcommand{\cclo}{cclo\xspace}
\newcommand{\qcclo}{$Q$-cclo\xspace}
\newcommand{\qcclos}{\qcclo's\xspace}
\newcommand{\smc}{smc-}

\newcommand{\csum}[2]{\osum{#1} + #2 + \sumo{#1}}

\newcommand{\cLpsum}{\csum Lp}

\newcommand{\li}[2]{{}_{#1}{#2}}
\newcommand{\csym}[2]{{\textstyle\sum}{#1}_i+ #2 +{\textstyle\sum^*}\li i{#1}}
\newcommand{\csymp}[1]{\csym{#1}{p}}

\newcommand{\cLpsym}{\csym Lp}

\newcommand{\vor}{\vartriangleleft}
\newcommand{\nvor}{\ntriangleleft}

\newcommand{\setof}[2]{\left\{ #1 \suchthat #2 \right\}}
\newcommand{\seqof}[2]{\left\< #1 \suchthat #2 \right\>}

\newcommand{\bbN}{\ensuremath{\mathbb{N}}\xspace}
\newcommand{\bbQ}{\ensuremath{\mathbb{Q}}\xspace}

\newcommand{\bbR}{\ensuremath{\mathbb{R}}\xspace}

\newcommand{\cC}{\mathcal{C}}

\newcommand{\cF}{\mathcal{F}}

\newcommand{\al}{\alpha}

\newcommand{\om}{\omega}
\newcommand{\si}{\sigma}

\newcommand\suchthat{\mathrel{:}}
\newcommand\val{v}
\newcommand\qe[1]{\exists #1}
\newcommand\qa[1]{\forall #1}
\newcommand\limp{\to}
\newcommand\subval{\mathrm{Val}}

\newcommand\vhi{\varphi}

\def\itm#1 {\item[(#1)]}
\def\<{\langle}
\def\>{\rangle}

\newcommand{\omG}{\omega_1^{\it G}}
\newcommand{\omGCB}{\omega_1^{\it GCB}}
\newcommand{\pfeil}[1]{\mathrel{\mathop{\longrightarrow}\limits^{#1} }}

\newcommand{\vecq}{\vec{q}\,}

\title{Continuous Fra{\"\i}ss{\'e} Conjecture}
\author{Arnold Beckmann \and Martin Goldstern \and Norbert Preining\\[2ex]
Institute of Discrete Mathematics and Geometry\\
  Vienna University of Technology, Austria\\
  \protect{\texttt{martin.goldstern@tuwien.ac.at}}\\
  \protect{\texttt{\{beckmann,preining\}@logic.at}}}
\begin{document}

\maketitle

\begin{abstract}
We will investigate the relation of countable closed linear orderings
with respect to continuous monotone embeddability and will show
that there are exactly $\aleph_1$ many equivalence classes with
respect to this embeddability relation. This is an extension of
Laver's result \cite{Laver:1971}, who considered (plain) embeddability,
which yields  coarser equivalence classes. Using this result we show that
there are only $\aleph_0$ many different Gödel logics.
\end{abstract}

\section{Introduction}

The starting point of the present work was the question `How many Gödel
logics are there?' 
This question led us to the study of embeddability relations of 
(countable) linear orderings. 
The most important result in this field is Laver's classical result 
on the Fra{\"\i}ss{\'e} Conjecture \cite{Laver:1971} which counts
the number of scattered linear orderings with respect
to bi-embeddability. 

We will generalize Laver's method to deal not only with monotone but
with continuous monotone embeddings, and come back to Gödel logics in
Section~\ref{sec:goedel}, where we use this result to compute the
number of Gödel logics. 
Gödel logics form a class of many-valued logics, which are one of
the three fundamental $t$-norm based logics.

% ***
% ***
% ***

% ***
% ***

% ***
% ***
% ***

Our main result is that the set of countable closed linear orderings 
is better-quasi-ordered by strictly monotone continuous embeddability, even
when we consider labeled countable closed linear orderings.
As a corollary we derive that there are only countably many Gödel logics.

% ***
% ***
% ***
% ***
% ***
% ***

% ***
% ***
% ***

The main concepts in all these discussions are `well-quasi orderings'
and `better-quasi-ordering', which have been introduced by
Nash-Williams in a series of five papers in the 1960s
\cite{Nash:1963,Nash:1964,Nash:1965a,Nash:1965b,Nash:1968} 

While considering embeddability relation of orderings, examples of
infinite descending sequences, as well as infinite antichains can be
given \cite{Dushnik,Sierpinski}. In \cite{Fraisse}, Fra{\"\i}ss{\'e}
made conjectures to the effect that the embeddability relation is more
well behaved in the case of countable order types (later extended to
scattered order types), stating that `every descending sequence of
countable order types is finite, and every antichain of countable
order types is finite.' This conjecture was finally proved by Laver
\cite{Laver:1971}.

\subsection{Basic concepts}

In our exposition we will mainly follow Rosenstein's textbook on linear
orderings \cite{Rosenstein:1982}, especially Chapter~10. To keep this
article self-contained we will give all the necessary definition and
cite some results, but ask the reader to consult the mentioned book
for motivation, background and history of these concepts and results, as well 
as for the proofs.  
% ***

\begin{definition} (\cite{Rosenstein:1982}, 10.12-10.15)
  A \emph{quasi-ordering} is a reflexive and transitive binary
  relation $\le_Q$ on a set $Q$. 
With $<_Q$ we denote the strict part of $\le_Q$, i.e.\ $p\mathrel{<_Q}q$
iff $p\mathrel{\le_Q} q$ and $q\mathrel{\nleq}_Q p$.
We will often drop the index ${}_Q$ if there is not danger of confusion.

% ***

We write $p \equiv_Q q$ iff both $p\mathrel{\le_Q} q$ and % ***
$q\mathrel{\le_Q} p$  hold.  
This is an equivalence relation; 
we write $Q/{\equiv}$ for the set of equivalence classes.

An infinite sequence  $\vec p  = \seqof{p_n}{n<\omega}$ is called good
if there are indices $n<k$ with $p_n\le p_k$; $\vec p $ is called 
bad if it is not good. 
$\vec p$ is called an 
 \emph{infinite descending chain} if
 $p_0\mathrel{>_Q}p_1\mathrel{>_Q}p_2\mathrel{>_Q}\dots$.
 It is called an \emph{anti-chain} of $Q$ if neither
 $p_i\mathrel{\le_Q} p_j$ nor $p_j\mathrel{\le_Q} p_i$ for $i\neq j$.

% ***
% ***

% ***
% ***
% ***

A set $Q$ is a \emph{well-quasi-ordering}, denoted \emph{wqo}, if
any/all conditions in Lemma~\ref{rosenstein:10.17} hold.
% ***
% ***
% ***
\end{definition}

\begin{lemma}\label{rosenstein:10.17} 
  {\normalfont (\cite{Rosenstein:1982}, 10.16--10.17)}
  Let $(Q,\le)$ be partial order. Then the following are equivalent: 
\begin{enumerate}
\item All sequences $ \vec q = \seqof{q_i}{i<\omega}$ are good. 
\item For all sequences $ \vec q = \seqof{q_n}{n<\omega}$ there
  is an infinite subsequence $\<q_n\suchthat n\in I\>$ which
  is either strictly increasing ($n<m$ implies $q_n<q_m$) or constant
  ($n<m$ implies % ***
   $q_n\equiv q_m$).
\item There are no infinite antichains and no infinite decreasing chains
 in $Q$. 
\end{enumerate}
\end{lemma}

\begin{definition} (\cite{Rosenstein:1982}, 10.19)
  Given quasi-orderings $Q_1$ and $Q_2$, we define the quasi-ordering
  $Q_1\times Q_2$ by stipulating that
  $\langle p_1,p_2\rangle\le\langle q_1,q_2\rangle$ if 
  $p_1\mathrel{\le_{Q_1}}q_1$ and $p_2\mathrel{\le_{Q_2}}q_2$.
\end{definition}

\begin{lemma} {\normalfont (\cite{Rosenstein:1982}, 10.20)}
  If $Q_1$ and $Q_2$ are wqo, then so is $Q_1\times Q_2$.
\end{lemma}

\begin{definition}\label{def:seq-order} (\cite{Rosenstein:1982}, 10.21, 10.24)
  Given a quasi-ordering $Q$, we define the quasi-ordering
  $Q^{<\omega}$, whose domain is the set of all finite sequences of
  elements of $Q$, by stipulating that
  $\<p_0,p_1,\dots,p_{n-1}\>\le\<q_0,q_1,\dots,q_{m-1}\>$ if there
  is a strictly increasing $h:n\to m$ such that
  $a_i\mathrel{\le_Q}b_{h(i)}$ for all $i<n$.

  We define the quasi-ordering $Q^{\omega}$ of $\omega$-sequences of
  elements of $Q$ by saying that
  $\seqof{p_n}{n<\omega}\le\seqof{q_n}{n<\omega}$ if there is a 
  strictly increasing $h: \omega\to\omega$ such that
  $a_n\mathrel{\le_Q}b_{h(n)}$ for all $n<\omega$.
\end{definition}

\begin{theorem}  {\normalfont (\cite{Rosenstein:1982}, 10.23)}
  If $Q$ is a wqo, then so is $Q^{<\omega}$.
\end{theorem}

\begin{definition}\label{def:barrier} (\cite{Rosenstein:1982}, 10.31--10.33)
If $c$ is a finite subset of $\bbN$, $d$ is any subset of $\bbN$, then we 
say that $d$ extends $c$ iff: $c = \{i\in d: i \le \max c\}$, i.e., 
if $c$ is an initial segment (not necessarily proper) of $d$.

  An infinite set $B$ of finite subsets of $\bbN$ is a \emph{block} if
  every infinite subset~$X$ of $\bigcup B:= \bigcup\setof{b}{b\in B}$ 
  has an initial segment in~$B$;
  that is, % ***
% ***
  $X$ extends some element in $B$.
  A block~$B$ 
% ***
  is called a \emph{barrier} if no two
  elements of~$B$ are comparable w.r.t.\ inclusion.

  A precedence relation $\vor$ on a barrier $B$ is defined as follows:
  if $b_1$ and $b_2$ are elements of $B$, then we say that $b_1$
  \emph{precedes}
  $b_2$, written $b_1\vor b_2$, if there are $i_1<i_2<\dots<i_m$ such
  that $b_1=\{i_1,i_2,\dots,i_k\}$ and $b_2=\{i_2,\dots,i_m\}$ for some
  $k$, $1\le k<m$.    
  (In particular, $\{i\} \vor \{j\}$ holds for all $i\not=j$.)

  A function $f:B\to Q$ on a barrier $B$ is \emph{bad} if, whenever 
  $b_1,b_2\in B$ and $b_1\vor b_2$, $f(b_1)\nleq_Q f(b_2)$.
  Otherwise we say that $f$ is good.
\end{definition}

\begin{definition} (\cite{Rosenstein:1982}, 10.30)
  We say that $Q$ is a \emph{better-quasi-ordering}, denoted
  \emph{bqo}, if every $f:B\to Q$ is good, for every barrier $B$ of
  finite subsets of $\bbN$.
% ***
% ***
% ***
\end{definition}

\begin{remark}  Every bqo is a wqo. 
\end{remark}
\begin{proof} 
Use the % ***
barrier $B = \{\{n\}: n\in \bbN\}$. 
\end{proof}
\begin{theorem}\label{bqoomega} {\normalfont (\cite{Rosenstein:1982}, 10.38)}
  If $Q$ is a bqo, then $Q^{<\omega}$ and $Q^\omega$ are bqo's.
\end{theorem}

\begin{theorem}\label{rosenstein:10.40} 
  {\normalfont (\cite{Rosenstein:1982}, 10.40)}
  Let $B$ be a barrier and suppose that $B=B_1\cup B_2$ is a partition
  of~$B$. Then there is a sub-barrier $C\subseteq B$ such that
  $C\subseteq B_1$ or $C\subseteq B_2$.
\end{theorem}

This ends the definitions and results we will need from
\cite{Rosenstein:1982}.

\begin{definition}
  A \emph{countable closed linear ordering}, denoted \emph{\cclo{}}, is
  a countable closed subset of $\bbR$.

  A \emph{strictly monotone continuous embedding} $h$
  (denoted \emph{\smc embedding}) from a \cclo $Q_1$ to a \cclo $Q_2$
  is an embedding $h: Q_1\to Q_2$ which is continuous on $Q_1$, i.e.\ 
  whenever $(p_n)_{n\in\bbN}$ is a sequence in $Q_1$ converging to an
  element $p$ in $Q_1$, then $(h(p_n))_{n\in\bbN}$ is a sequence in $Q_2$
  converging to an element $h(p)$ in $Q_2$,
  and strictly monotone on $Q_1$, i.e.\ whenever $p,q\in Q_1$ with 
  $p\mathrel{<_{Q_1}}q$ then $h(p)\mathrel{<_{Q_2}}h(q)$. 
  (Here, ``convergence'' is always understood as convergence in 
   the usual topology of $\bbR$.)
\end{definition}

\begin{definition}[labeled \cclo]
  In addition to \cclo, we will also have to consider the following
  notion:  Fix a quasi-order $Q$ (usually a bqo, often  a
  finite set or an ordinal). A \emph{\qcclo} is a function $A$ whose domain
  $\dom{A}$ is a \cclo\ and whose range is contained in $Q$.

  We write $A \contemb B$ ($A$ is $Q$-\smc embeddable into $B$, or shortly
  $A$ is embeddable into $B$) iff there is a \smc embedding
  $h$ from $\dom{A}$ to $\dom{B}$ with the property 
  $A(a) \le_Q B(h(a))$ for all $a\in\dom{A}$.
\end{definition}

% ***
% ***
% ***
% ***

If $Q$ is a singleton, then $A\contemb B$ reduces just to a  
\smc embedding from $\dom{A}$ to $\dom{B}$.  
If $Q = \{p,q\}$ is an antichain, or satisfies $p<q$, 
and $A(0)=A(1)=q=B(0)=B(1)$, $B(b)=p$ for all $b\not=0,1$, then
$A\contemb B$ means that there is a \smc embedding
from $\dom{A}$ to $\dom{B}$ which moreover preserves $0$ and $1$.
Such embeddings will play an important r\^ole when we investigate Gödel sets
and the number of Gödel logics.

% ***
\section{$Q$-labeled countable closed linear orderings}

Let us fix some bqo $(Q,\le)$ for defining \qcclos.

\begin{notation}
  We will use the following notation throughout the paper:
  \[ \cLpsum \]
  or 
  \[ \cLpsym \]
  When we write this term the following conditions are imposed:
  \begin{itemize}
  \item $p$ is an element of $Q$.
  \item All the $L_i$ and $"iL$ are \qcclos. 
  \item Either all $L_i$ are empty, or none of them are empty.
      Similarly, either all $"iL$ are empty, or none of them are.  
      We do not allow all $L_i$ {\em and} all $"iL$ to be empty. 
   \item $\dom L_i<\dom L_{i+1} < \dom "{i+1}L<\dom "iL$ for all $i$,
    where we write $A < B $ for ``$A=\emptyset \vee B=\emptyset \vee \sup A < \inf B$''. 
    % ***
    % ***
In particular, between the domains of any two of them (in the non-empty case)
we can find an open interval.
  \item $\lim_{n\to\infty} a_n = \lim_{n\to\infty} "na$, 
    whenever  $a_n\in \dom L_n$ and $"na\in\dom "nL$.
  \end{itemize}
  % ***
  % ***
  % ***
  % ***
  % ***

The {\em meaning} of such a term is the $Q$-cclo  $L$ whose domain is the 
set $\bigcup_i L_i \cup \{x\} \cup \bigcup_i "iL$ (where
$x= \lim_{n\to\infty} a_n $ and/or  $x  = \lim_{n\to\infty} "na$ for any/all 
sequences satisfying  $a_n\in \dom L_n$ and $"na\in\dom "nL$), and 
the function
$L$ extends all functions $L_i$ and $"iL$, and $L(x) = p$.

A ``finite sum'' $$L = L_1 + \cdots + L_n$$ is defined naturally:  we allow 
this expression only when all $L_i$ are nonempty and satisfy $\max\dom L_i< 
\min \dom L_{i+1}$.   In this case we let
$\dom( L) = \bigcup_i\dom( L_i)$ and $L = \bigcup_i L_i$. 
\end{notation}

We will consider two slightly different operations ($S$, $S'$ below)
to build complicated \qcclos from simpler ones.   
These two operations naturally correspond to two notions $\rk$, $\rk'$ of rank;
a third rank that we occasionally use is the classical Cantor-Bendixson rank 
$\cbrk$ of a \cclo.

\begin{definition}
  Let $\OO$ be a class of \qcclos. We let $S(\OO)$ 
  (`sums from $\OO\/$\/') be the set of all \qcclos which are
finite sums of \qcclos from $\OO$, plus the set of all \qcclos
 of the form 
  \[  \cLpsum \]
  where $p\in Q$ and all $L_n$ and all $"nL$ are in $\OO$.

  We let $S'(\OO)$ (`unbounded sums from $\OO$') be the set
  of all \qcclos of the form 
  \[  \cLpsum \]
  where $p\in Q$ and all $L_n$ and all $"nL$ are in $\OO$, and 
  \[ \forall n \, \exists k>n \ L_n\contemb L_k 
  \mbox{ \  and  \ } 
  \forall n \, \exists k>n \ "nL\contemb "kL. \]
\end{definition}

As a consequence of the above definition we obtain for unbounded sums, 
that for all~$n$ there are infinitely many $k>n$ such that $L_n\contemb L_k$ 
and $"nL\contemb "kL$.

% ***
% ***
% ***
% ***
% ***

% ***
% ***
% ***

% ***
% ***
% ***
% ***

% ***
% ***
% ***
% ***

\begin{definition}\label{def:abcd} \ 
\begin{enumerate}
\itm a 
  Let $\sC$ be the set of all 
% ***
% ***
\qcclos.

\itm b 
  Let $\sC_0= \sC_0' $ be the set of all \qcclos with singleton or empty 
  domain.
  For any $\alpha \le \omega_1$ let 
 \[ \sC_{\alpha+1} = S(\sC_\alpha) \cup \sC_\alpha \qquad 
  \sC'_{\alpha+1} = S'(\sC'_\alpha) \cup \sC'_\alpha \]
and for limit ordinals $\delta>0 $ let  
$\sC_\delta = \bigcup_{\alpha<\delta} \sC_\alpha$, 
$\sC'_\delta = \bigcup_{\alpha<\delta} \sC'_\alpha$.

\itm c 
  For any % ***
  $L\in\bigcup_\al\sC_\al$ 
  we define the \emph{rank} of~$L$ ($\rk(L)$) as
  the first ordinal~$\al$ at which $L$ occurs in $\sC_{\al+1}$.
  Similar, we define $\rk'(L)$ for $L\in\bigcup_\al\sC'_\al$ as
  the first ordinal~$\al$ at which $L$ occurs in $\sC'_{\al+1}$.
% ***

\itm d 
  The set of all \qcclos whose domains are suborderings of $\dom L$
  is  denoted with $\sC(L)$.
\end{enumerate}
\end{definition}

It is clear that $\sC'_{\omega_1} \subseteq \sC_{\omega_1} \subseteq \sC$. 
We will show that $\sC = \sC_{\omega_1}$, and that every order in $\sC$ can 
be written as a finite sum of orders from $\sC'_{\omega_1}$.

\begin{lemma}
 % ***
  $ \sC = \sC_{\omega_1}$.   
That is, for every \qcclo\ $L$ there is a countable
ordinal $\alpha$  such that $L\in \sC_{\alpha}$.
% ***
\end{lemma}

\begin{proof}
We use the Cantor-Bendixson decomposition, more precisely we 
use induction on the Cantor-Bendixson rank of $V=\dom L$. 

For every scattered closed set $V$  
there is an ordinal $\cbrk(V)$ (the Cantor-Bendixson rank
of $V$) and a decomposition 
\[ V = \bigcup _{\alpha\le \cbrk(V)} \CB_\alpha(V),\]
where $\CB_0(V) $ is the set of isolated points of $V$, and
more generally each set $\CB_\alpha(V)$ is the set of isolated points
of   $V\setminus \bigcup_{\beta < \alpha} \CB_\beta(V)$, and 
$\CB_{\cbrk(V)}(V)$ is finite and nonempty. 

Assume for the moment that  $\CB_{\cbrk(V)}(V)$  
is a singleton $\{x^*\}$. 
If $\cbrk(V) = 0$, then % ***
$L\in \sC_0$. 
If $\cbrk(V) > 0$, fix an increasing sequence $\<x_n\>$ and 
a decreasing sequence $\<"nx\>$, both with limit $x^*$,
and $x_n, "nx\notin V$.  Now it is easy to see that for 
all $\beta<\cbrk(V)$
\[ \CB_\beta( V \cap [x_n, x_{n+1}] ) = 
 \CB_\beta( V)  \cap [x_n, x_{n+1}],\]
so $\cbrk(V\cap [x_n,x_{n+1}]) < \cbrk(V)$, similarly for
 $V\cap ["nx,"{n+1}x]$.   Now we can 
use the induction hypothesis. 

If $\CB_{\cbrk(V)}(V)$  is not a singleton then we can  write $V = V_1 + 
\cdots + V_n$ for some finite $n$, with each $\CB_{\cbrk(V)}(V_k)$ a singleton,
then proceed as above. 
% ***
% ***
% ***
\end{proof}

% ***
% ***
% ***
% ***
% ***
% ***
% ***
% ***
% ***
% ***
% ***
% ***
% ***
% ***

% ***
% ***
% ***

\begin{definition}  
The set $\sC':= \sC'_{\omega_1}$ is the smallest family of \qcclos
which contains all the singletons and is closed under unbounded sums $S'$. 
% ***
% ***
\end{definition}

% ***

\begin{theorem}\label{thm:finite-ha}
  Let $L$ be a \qcclo and assume that $(\sC(L), \contemb)$ is a wqo. (See
  Definition~\ref{def:abcd}(d).)
  Then $L$ is a finite sum of % ***
elements in $\sC'$.
\end{theorem}

\begin{proof}
  Induction on $\rk(L)$: Assume that
  \[ L = \cLpsum \]
  where all the $L_i$ and $"iL$ are in $\sC'$.
  Suppose that, for all but a finite number of $L_i$, each
  $L_i$ is embeddable in infinitely many $L_j$, and for all but a
  finite number of $"iL$, each $"iL$ is embeddable in infinitely
  many $"jL$. Then $L$ can be written as
  \[ L_0 + \dots L_{k-1} + (L_{k+0} + L_{k+1}
  +\dots+p+\dots+"{l+1}L+"{l+0}L) + "{l-1}L+ \dots +"0L \] where
  each summand is in $\sC'$.

  Otherwise there are either infinitely many $L_i$ or $"iL$ each
  embeddable in only finitely many $L_j$ or $"jL$, resp. We then
  find a either a subsequence $\<L_{h(n)}\suchthat n<\omega\>$ or
  $\<"{h(n)}L\suchthat n<\omega\>$ no entry of which can be embedded
  in any subsequent entry. This bad sequence of suborderings of $L$
  contradicts the hypothesis of the theorem.
\end{proof}

\begin{theorem}\label{thm:bqo-wqo}
  If $(\H,\contemb)$ is a bqo, then $(\sC,\contemb)$ is a wqo.
\end{theorem}

\begin{proof}
  We will show for all countable $L$ by induction on the 
    rank~$\rk(L)$ (that is the rank w.r.t.\ the classes in $\sC$ as defined
    in Definition \ref{def:abcd} (c)),
  that the collection $\sC(L)$ of \qcclos whose domains are suborderings
  of $\dom L$ is a wqo w.r.t \smc embeddability.

  First we show that, if $K$ is in $\sC(L)$, then $K$ can be
  written as $K=\csymp J$, where all the $J_i$ and $"iJ$ are
  in $\sC'$. To prove this, observe that $L$ can be written
  as $\cLpsym$ where the ranks of the $L_i$ and $"iL$ are strictly less
  than the rank of $L$. Using the induction hypothesis, we see that
  $(\sC(L_i),\contemb)$ and $(\sC("iL),\contemb)$ are wqo. If $\dom K$ is a
  sub ordering of $\dom L$, it can be written as $K=\csym Kq$ 
% ***
% ***
% ***
  with $K_i\in\sC(L_i)$ and $"iK\in\sC("iL)$.
  Thus, by Theorem~\ref{thm:finite-ha}, each $K_i$ and $"iK$
  can be written as finite sum of elements $J_j$ and
  $"jJ$ in $\sC'$, and $K$ as $\csum Jq$.

% ***
% ***

  Now consider  a sequence $\<K^l\suchthat l<\omega\>$, where each $K^l$
  is a \qcclo and subordering of $L$.  We will repeatedly thin out
  this sequence, eventually arriving at a sequence which is good, which 
  will show that our original  sequence was good.
  After having thinned
  out the sequence $\<K^l\suchthat l<\omega\>$ to a sequence 
  $\<K^{l_i}\suchthat i<\omega\>$, we will (for notational simplicity)
  relabel our index set so that we will also call the new sequence 
  $\<K^l\suchthat l<\omega\>$. 

  Each $K^l$ can be written as
  \[ K^l = \csum {J^l}{p^l} \]
  where each of the summands is in $\sC'$.
  Using Lemma~\ref{rosenstein:10.17} % ***
  we thin out our sequence to a new sequence 
  (again called $\<K^l\suchthat l<\omega\>$) such that 
  $p^j \mathrel{\le_Q} p^k$ for all $j<k$. 

  By Theorem~\ref{bqoomega} we know that $\H^\omega$ is a bqo, in
  particular a wqo.
% ***
% ***
% ***
  Consider the $\omega$-tuples
   ${\cC}^l = \< J_0^{l}, J_1^{l}, \dots \>\in \H^\omega$.
  Using Lemma~\ref{rosenstein:10.17} % ***
  we  can thin out our sequence to obtain a sequence satisfying 
  ${\cC}^{j} \contemb {\cC}^{k}$  for any $j<k$.

  We now apply the fact that $\H^\omega$ is wqo to the sequence  
  ${}^n\cC = \< "0J^n, "1J^n, \dots \>\in \H^ \omega$ to see 
  that without loss of generality we may also assume 
  ${}^j{\cC} \contemb {}^k{\cC} $ for all $j<k$. 

  Now pick any $n<m$, and consider the sums
  \[ K^n = \csum {J^n}{p^n} \]
 and 
  \[ K^m = \csum {J^m}{p^m} .  \]
 Write $x^n$ and $x^m$ for the central points of $K^n$ and $K^m$, respectively
 (i.e., $x^n = \sup_n \bigcup_i \dom J^n_i = \inf_n \bigcup_i \dom "i{J^n}$, etc.)

 We know $p^n\mathrel{\le_Q} p^m$,  $\cC^n \contemb \cC^m$, 
                                    $"n\cC \contemb "m\cC$.

  Thus, there are strictly increasing functions $g$ and $h$ from
  $\bbN$ to $\bbN$, 
  such that for all $i$,  
  $J_i^n \contemb J_{g(i)}^m$ and $"i{J^n} \contemb "{h(i)}{J^m}$.  Let $\alpha_i$ and $"i\alpha$ be functions that witness this, i.e., 
  let $\alpha_i$ be a function mapping $\dom J_i^n$ to $\dom J_{g(i)}^m$ 
with $J_i^n(x) \le J_{g(i)}^m(\alpha_i(x))$ for all $x\in \dom J_i^n$, 
and similarly  $"iJ^n(x) \le "{g(i)}J^m("i\alpha(x))$ for all $x\in \dom "iJ^n$, 

Now define $\alpha:\dom K^n \to \dom K^m$ naturally:  $\alpha$ extends all 
functions $\alpha_i$ and $"i\alpha$, and $\alpha(x^n) = x^m$.  Clearly
$\alpha$ witnesses $K^n\contemb K^m$. 

\medskip

  Finally, if $\setof{K_i}{i<\omega}$ is an arbitrary sequence,
  where each $K_i$ is in $\sC$, then each $K_i\in\cC(K)$ where 
  $K= \csum {K}{p}$ for arbitrary $p$ and empty $"iK$. According to
  the above remarks, the sequence $\setof{K_i}{i<\omega}$ must be
  good, so that $\cC$ is a wqo.
\end{proof}

% ***
% ***
% ***

% ***
% ***
% ***
% ***
% ***
% ***
% ***

\begin{theorem}\label{thm:ha-bqo}
  $(\H,\contemb)$ is a bqo.
\end{theorem}

We prove the Theorem by a series of lemmas.
The first lemma holds for general quasi-orderings which are equipped
with a rank function, it forms the main technical part of the proof 
of Theorem~\ref{thm:ha-bqo}.

% ***
% ***
% ***
% ***
% ***
% ***
% ***
% ***
% ***
% ***
% ***
% ***
% ***
% ***
% ***
% ***

Let $(Q,\le)$ be a quasi-ordering, 
and let $\rho$ be a rank function from $Q$ into the ordinals
(i.e., a function satisfying $\rho(x)\le \rho(y)$ whenever 
$x\le y$). 
Let $\cF$ denote the set of all functions $g:B\to Q$
where $B$ is a barrier of finite subsets of $\bbN$. (See 
Definition~\ref{def:barrier}.)

We say that $C$ \emph{is an extended sub-barrier of} $B$
if $\bigcup C\subseteq\bigcup B$ and 
if every element of $C$ is an extension (not necessarily proper)
of an element of $B$.
$C$ is called a \emph{proper} extended sub-barrier of $B$
if $C$ {is an extended sub-barrier of} $B$ and 
at least one element of $C$ properly extends some element of $B$.
For two functions $g:B\to Q$ and $h:C\to Q$ in $\cF$
we say that $h$ \emph{is shorter than} $g$ if
$C$ is a proper extended sub-barrier of $B$ and
if $g$ and $h$ coincide on $B\cap C$, and 
if, whenever $c\in C$ properly extends $b\in B$, 
$h(c)\le g(b)$ and $h(c)$ has lower rank than $g(b)$.
The following Lemma can be extracted from the proof of Theorem 10.47
in Rosenstein \cite{Rosenstein:1982}. 
Recall from Definition~\ref{def:barrier} that a function $f: B\to Q$
is called \emph{bad}
 if, whenever $b_1,b_2\in B$ and $b_1\vor b_2$,
$f(b_1)\nleq f(b_2)$.

\begin{lemma}\label{lem:minimal-bad}
  If $\cF$ contains some bad function, than it contains 
  some minimal bad function,
  i.e.\ one which is minimal w.r.t.\ `shorter'.
\end{lemma}

\begin{proof}
Assume for the sake of contradiction that 
$\cF$ contains some bad function, 
but for any bad $g\in\cF$ there is some bad $h\in\cF$ 
which is shorter than $g$.

Let $g:B\to Q$ be bad.
With $k(g)$ we denote the minimal $k$ such that there is a shorter
$h:C\to Q$ and a $b\in B$ which is properly extended by some element in $C$
with $\max b\le k$.
Fix some witnesses $C$, $h$ and $b$ for $k(g)$.
We define $D$ as the set of all $d\in B$ which do not have extensions
in $C$ and which fulfill $d\subset[0,k(g)]\cup\bigcup C$.
Obviously $C\cap D=\emptyset$.

% ***
% ***

First observe that for $d\in D$ we have $d\not\subset\bigcup C$:
Assume for the sake of contradiction that $d\subset\bigcup C$.
Let $X$ be the infinite set $d\cup\big(\bigcup C\cap[\max d,\infty)\big)$,
then $X\subseteq\bigcup C$, hence there is some $c\in C$ which is 
extended by $X$.
Since $X$ is also  an extension of $d$,  $c$ extends $d$ or vice versa.
As $c$ extends some element in $B$ and $d\in B$, 
we have that $c$ cannot be properly extended by $d$ because $B$ is a barrier.
But by definition of $D$ we also have that $c$ does not extend $d$.
Contradiction.

% ***

Now, $B^*:=C\cup D$ is a barrier and $g^*:B^*\to Q$ defined by
$g^*(c)=h(c)$ for $c\in C$ and $g^*(d)=g(d)$ for $d\in D$
is bad and shorter than $g$.   

We verify these claims:
First note that $\bigcup B^* \subseteq \bigcup C\cup [0,k(g)]$. 
For $B^*$ to be a block let $X\subseteq\bigcup B^*$ be infinite.
There is some $d\in B$ which is extended by $X$ (as $B$ is a 
block, and $\bigcup B^* \subseteq \bigcup B$).
If $d$ is not already in $B^*$ then, 
by definition of $D$, $d$ has some extension 
in $C$ which must be proper as $d\notin C$.
Thus $d\subset\bigcup C$ and $\max d\ge k(g)$, 
hence $X\subseteq\bigcup C$ % ***
as $\bigcup B^*\subseteq\bigcup C\cup [0,k(g)]$.
But then there is some $c\in C$ which is extended by $X$.
Altogether this shows that $B^*$ is a block.

Assume that $B^*$ is not a barrier, then there must be $c\in C$ and
$d\in D$ which are comparable.
As $c$ is the extension of some element in $B$ and $d\in B$, 
we have $c\nsubseteq d$ because $B$ is a barrier.
But $d\in D$ implies $d\not\subset\bigcup C$, hence $d\nsubseteq c$.
Contradiction.
Hence $B^*$ must be a barrier.

Obviously, $g^*$ is shorter than $g$, as $h$ already has been shorter 
than $g$.
To verify that $g^*$ is bad we assume for the sake of contradiction
that $c_1\vor c_2$ and $g^*(c_1)\le g^*(c_2)$.
As $h$ is bad, $c_1$ and $c_2$ cannot be in $C$ at the same time.
Similar with $g$, they cannot be in $D$ at the same time.
If $c_1\in C$ and $c_2\in D$, 
then $c_2\not\subset\bigcup C$
which together with $c_1\vor c_2$, $c_1\in C$ and the definition of $D$
shows $\max c_1<k(g)$, hence $c_1\in B$.
Hence $g(c_1)=h(c_1)=g^*(c_1)\le g^*(c_2)=g(c_2)$ 
contradicting that $g$ is bad.
Therefore, $c_1\in D$ and $c_2\in C$.
There is some $b_2\in B$ such that $b_2$ is extended by $c_2$.
If $c_1\nvor b_2$ then $b_2\subsetneq c_1$ 
which contradicts that $B$ is a barrier.
Hence we have $c_1\vor b_2$.
But then $g(c_1)=g^*(c_1)\le g^*(c_2)=h(c_2)\le g(b_2)$
contradicts that $g$ is bad.
Altogether this shows that $g^*$ is bad.

We now define a sequence of bad elements $f_n\in\cF$ in the following way.
Let $f_0:B_0\to Q$ be some bad element in $\cF$, and define recursively
$B_{n+1}:=B_n^*$ and $f_{n+1}:=f_n^*$.
Let $k_n:=k(f_n)$.
Then $k_{n+1}\ge k_n$ because `shorter' is transitive and $k_n$ is chosen
minimal.
Furthermore, $k_n=k_m$ for only finitely many $m$ since
$\setof{b\in B_n}{\max b=k_n}$ is finite.
Hence $\<k_n \suchthat n<\om\>$ is a non-decreasing unbounded sequence of
natural numbers.  
Also observe that if $b\in B_n$ and $\max b<k_n$ and $n<m$ then $b\in B_m$,
and if $b\in B_m\cap B_n$ then $f_m(b)=f_n(b)$.

Let $B:=\bigcup\setof{\bigcap\setof{B_n}{n\ge m}}{m<\om}$.
We show that $B$ is a barrier.
Let $M:=\bigcap\setof{\bigcup B_n}{n<\om}$.
$M$ is infinite because $k_n\in M$ for all $n$.
Let $X\subseteq M$ be infinite.
Then for all $n<\om$ we have $X\subseteq\bigcup B_n$, hence there is
some $b_n\in B_n$ which is extended by $X$.
If $b_{n+1}$ is a proper extension of $b_n$ then the rank of $f_{n+1}(b_{n+1})$
is strictly smaller than the rank of $f_n(b_n)$, hence, for some $m$,
$b_n=b_m$ for all $n\ge m$,
i.e.\ $b_m\in\bigcap\setof{B_n}{n\ge m}\subseteq B$.
In particular, $M\subseteq\bigcup B$ by taking $X:= M\cap[m,\infty)$ 
for $m\in M$.
If $k\in\bigcup B$, then there is some $b\in B$ with $k\in b$.
$b\in B$ implies that there is some $m$ with $b\in\bigcap_{n\ge m}B_n$.
Thus $k\in\bigcup B_n$ for all $n\ge m$. 
Also
$k\in\bigcup B_m \subseteq \bigcup B_{m-1} \subseteq\dots\subseteq\bigcup B_0$,
hence $k\in M$.
This shows $\bigcup B\subseteq M$.
Thus $M=\bigcup B$ and $B$ is a block.
Let $b,c\in B$, then $b,c\in B_n$ for some $n$, hence they are not comparable
as $B_n$ is a barrier.
Altogether this shows that $B$ is a barrier.

For $b\in B$ let $m_b := \min\setof{m}{b\in\bigcap\setof{B_n}{n\ge m}}$.
We define $f:B\to Q$ by $f(b):=f_{m_b}(b)$
and show that $f$ is minimal w.r.t.\ `shorter' and bad.
$f$ is shorter than $f_n$ for all $n$,
because `shorter' is transitive,
$B$ is an extended sub-barrier of $B_n$,
if $b\in B\cap B_n$ then $m_b\le n$ hence $f(b)=f_{m_b}(b)=f_n(b)$,
and if $c\in B$ properly extends $b\in B_n$, then
$m_c>n$ and $f(c)=f_{m_c}(c)\le f_n(b)$ and 
$\rho(f(c))=\rho(f_{m_c}(c))<\rho(f_n(b))$.
$f$ is bad,
because if $b,c\in B$, and w.l.o.g.\ $m_b\le m_c$, then
$f(b)=f_{m_c}(b)$ and $f(c)=f_{m_c}(c)$ and $f_{m_c}$ is bad.
By our general assumption 
there is some bad $f':B'\to Q$ which is shorter than $f$.
Then there are $b'\in B'$ and $b\in B$ such that $b'$ properly extends $b$.
Choose $n$ with $k_n>\max b'$.
Now $f'$ is shorter than $f_n$ because $f$ is shorter than $f_n$ 
and `shorter' is transitive.
But this contradicts the minimality of $k(f_n)$.
Hence our general assumption has been wrong, and the theorem is proved.
\end{proof}

Recall that the rank $\rk'(L)$ of $L\in\sC'$ is given 
by the minimal $\al$ such that $L\in\sC'_{\al+1}$.
A \emph{$\sC'$-term} for $L\in\sC'$ with $\rk'(L)>0$ is a faithful witness 
for $L\in\sC'$, i.e.\ a decomposition $L=\cLpsym$ with
all the $L_i$ and $\li iL$ in $\sC'$ and 
$\rk'(L_i)<\rk'(L)$ and $\rk'(\li iL)<\rk'(L)$ for all $i$.

\begin{lemma}\label{lem:charc-indec}
  Let $L$ and $K$ be in $\sC'$, with
  $L=\cLpsym$ and $K=\csym Kq$ being $\sC'$-terms of them.
  If $p\le q$ and each $L_i$ is  embeddable into some $K_j$ and
  each $\li iL$ is embeddable into some $\li jK$, then $L\contemb K$.
\end{lemma}

\begin{proof}
Let the assumptions of the lemma be fulfilled.
Then there are
$k_i,l_i\in\om$ such that $k_i<k_{i+1}$, $l_i<l_{i+1}$ and
$L_i\contemb K_{k_i}$ and $"iL\contemb "{l_i}K$
because $L,K$ are in $\sC'$.
Fix \smc embeddings $\si_i:\dom L_i\to \dom K_{k_i}$ and 
$"i\si:\dom "iL\to \dom "{l_i}K$ 
witnessing $L_i\contemb K_{k_i}$ resp.\ $"iL\contemb "{l_i}K$,
and let $b := \lim_i(\sup \dom L_i)$ and $c := \lim_i(\sup \dom K_i)$.
We define a map $\si:\dom L\to \dom K$ by
\[  
  \si(a) := 
  \begin{cases}
    \si_i(a) & \text{if $a\in \dom L_i$}\\
    "i\si(a) & \text{if $a\in \dom "iL$}\\
    c        & \text{if $a=b$}
  \end{cases}
\]
% ***
% ***
% ***
% ***
% ***
% ***
Then $\si$ is a \smc embedding witnessing  $L\contemb K$.
% ***
% ***
% ***
% ***
% ***
% ***
% ***
% ***
% ***
% ***
% ***
% ***
% ***
% ***
% ***
% ***
% ***
% ***
% ***
% ***
\end{proof}

\begin{proof}[Proof of Theorem \ref{thm:ha-bqo}]
Assume for the sake of contradiction that $(\H,\contemb)$ is not a bqo.
By applying Lemma~\ref{lem:minimal-bad} we can find some
$f:B\to\sC'$ which is bad and minimal w.r.t.\ `shorter'.
For each $b\in B$ we fix some $\sC'$-term $f(b)=\cLpsym$.

For any  $a,b\in B$ with $a\vor b$ we have that 
$f(a)=\cLpsym\ \ncontemb\ \csym Kq=f(b)$,
hence, by applying Lemma~\ref{lem:charc-indec}, we see that at least one 
of the following holds: 
\begin{enumerate}
\itm i 
$p\nleq q$
\itm ii  for  some $i$: $L_i\ncontemb K_j$ for all $j$, 
\itm iii  for some $i$: $"iL\ncontemb "jK$ for all $j$.
\end{enumerate}
By applying Theorem~\ref{rosenstein:10.40} 
% ***
we can find
a sub-barrier $B'$ such that one the cases (i), (ii), (iii) 
 always happens on $B'$.
In the first case this would form a bad sequence in $(Q,\le)$ which would
contradict that $(Q,\le)$ is a wqo.
Thus, w.l.o.g.\ we may assume that for all $a,b\in B'$ with $a\vor b$
there is some $i$ such that $L_i\ncontemb K_j$ for all $j$.
Let $B'(2):=\setof{b_1\cup b_2}{b_1,b_2\in B'\text{ and }b_1\vor b_2}$,
then $B'(2)$ is an extended sub-barrier of $B$.
Define $g:B'(2)\to\sC'$ by letting $g(b_1\cup b_2)$ be the first $L_i$
in $f(b_1)=\cLpsym$ which is not embeddable into 
any $K_j$ from $f(b_2)=\csym Kq$.
Then obviously $g$ is shorter than $f$.
But also $g$ is bad, because if $b_1\cup b_2\vor b_3\cup b_4$ then $b_2=b_3$
and hence $g(b_1\cup b_2)\ncontemb g(b_3\cup b_4)$.
This contradicts the minimality of~$f$.
\end{proof}

Theorems \ref{thm:bqo-wqo} and \ref{thm:ha-bqo} together yield the
following result:

\begin{corollary} \label{cor:sc-wqo}
$(\sC,\contemb)$ is a wqo.
\qed
\end{corollary}

For the next corollary, we need the following two well-known properties of 
wqo's:

\begin{lemma}\label{omega1.wqo}
 Let $(Q, {\le})$ be a wqo
 with uncountable many $\equiv$-equivalence classes.
% ***
% ***
% ***
 Then there exists a  1-1 monotone map $f:\omega_1\to Q$. 
\end{lemma}
\begin{proof}[Proof sketch]
W.l.o.g.\ let each equivalence class of $Q/\equiv$ consist of one element.
If each uncountable subset $Q'\subseteq Q$ contains some
element $q$ such that also $\setof{r\in Q'}{q \nleq r}$ is uncountable,
then we can find sequence 
\[Q = Q_0  \supseteq Q_1\supseteq Q_2 \supseteq\cdots  \]
of uncountable sets with 
elements $q_n\in Q_n$, $Q_{n+1}:= \setof{r\in Q_n}{q_n \nleq r}$.  But then 
$q_n\nleq q_k$ for all $n<k$, contradicting  the assumption that $Q$ is
wqo.

So there must be an uncountable subset $Q'\subseteq Q$ such that
 for any $q\in Q'$, the set $\setof{r\in Q'}{q \nleq r}$ is countable.
But then we can easily find a copy of $\omega_1$ in $Q'$.
% ***

Alternatively, start with any 1-1 sequence $\<q_i:i\in \omega_1\>$
in $Q$; define
a coloring $f:[\omega_1]^2\to 2$ by $f(i<j) = 0$ iff $q_i < q_j$, and 
apply the Erd\hungarian os-Dushnik-Miller theorem $\omega_1\to (\omega_1,\omega)$.  (See \cite[Theorem 11.1]{EHMR}.)

% ***

% ***
% ***
% ***
% ***
% ***
% ***
% ***
% ***
% ***
% ***
% ***
% ***
% ***
% ***
% ***
% ***

% ***
% ***
% ***
% ***
% ***
% ***
% ***
% ***
% ***
% ***
% ***
% ***
% ***
% ***
% ***
% ***
% ***
\end{proof}

\begin{lemma}\label{lemma:Qomega-countable}
Let $Q$ be a countable bqo (or at least assume that $Q$
has only countably many $\equiv$-equivalence classes). 

Then $Q^\omega$ (quasiordered as in Definition~\ref{def:seq-order})
   has only countably many equivalence classes. 
\end{lemma}

\begin{proof} Part I:
 We first consider the set $Q^*$ of all sequences 
$\vecq  = \<q_0,q_1,\ldots\>\in Q^\omega$ satisfying 
\[ \forall k\,\exists n>k : q_k \le q_n.\]
and show that this set is countable (modulo $\equiv$). 

By Theorem~\ref{bqoomega},
$Q^\omega$ and hence also $Q^*$ is a wqo.  
  Assume that $Q^*$ has uncountably many 
$\equiv$-classes, then  by Lemma~\ref{omega1.wqo} we can find a sequence 
$\seqof{\vecq^i}{i \in \omega_1}$, 
\[ \vecq^i = \<q^i_0,q^i_1,\ldots\>\in Q^* \]
with $i< j \Rightarrow \vecq^i \le \vecq^j$,
$\vecq^j \nleq \vecq^i$.

Let $\al < \omega_1$ be so large such that every element of $Q$ which 
appears somewhere as $q^j_n$ is $\le$ to some $q^{j'}_{n'}$ with $j'<\al$.  

We claim that $\vecq^{\alpha+1} \le \vecq^\al$, which will
be the desired  contradiction.

By definition of $\alpha$, $\forall n \, \exists i < \al\, \exists n':  
 q^{\al+1}_n \le q^i_{n'}$.  So for every $n$ there is $n''$ with 
$q^{\al+1}_n \le q^\al_{n''}$. 
Using  $\vecq^\al\in Q^*$, we can find a sequence $k_0 < k_1 < \cdots $
with $q^{\al+1}_n \le q^\al_{k_n}$ for all $n$, which means 
 $\vecq^{\alpha+1} \le \vecq^\al$. 
\medskip

\noindent Part II:  For any sequence 
$\vecq = \<q_0,q_1,\ldots\>\in Q^\omega$
we can find a natural number $N=N_{\vecq}$ such that 
$\forall k\ge N \,\exists n>k : q_k \le q_n$, otherwise we get
(as in the proof of Theorem~\ref{thm:finite-ha}) a contradiction 
to our assumption that $Q$ is a wqo. 

Now assume that $Q^\omega / {\equiv}$ is uncountable, then we can 
find a natural number $N^*$ and  an
uncountable family $\seqof{\vecq^i}{i < \omega_1}$ of pairwise
nonequivalent sequences in $Q^\omega$ such that for all $i$, 
$N_{\vecq^i} = N^*$.   Moreover, we may assume that 
all initial segments  $\<q^i_0, \ldots, q^i_{N^*}\>$ are equal to
each other.     Consider the tails
$\<q^i_{N^*+1}, q^i_{N^*+2} , \ldots\> \in Q^\omega$.  By definition
of $N_{\vecq^i}$, these tails are all in $Q^*$, defined in part I, above. 

Hence we can find $i\not= j$ such that 
\[ \<q^i_{N^*+1}, q^i_{N^*+2} , \ldots\> \equiv 
\<q^j_{N^*+1}, q^j_{N^*+2} , \ldots\>.\]
But then also $\vecq^i \equiv \vecq^j$. 
\end{proof}

\begin{corollary}\label{cor:S-countable}
Assume that our basic wqo $Q$ is countable. 
Then, for any set  $\OO \subseteq \sC$   with $\OO/{\equiv}$ countable
we also have that  $S'(\OO)/{\equiv}$ and even 
$S(\OO)/{\equiv}$  are countable. 
\end{corollary}

\begin{proof}

 If $\<L_0,L_1,\ldots\> \contemb \<L'_0,L'_1,\ldots\>$
and 
 $\<"0L,"1L,\ldots\> \contemb \<"0L',"1L',\ldots\>$ and $p \le p'$, then
also 
$$ L_0 + L_1 + \cdots + p + \cdots "1L + "0L 
 \ \contemb\ 
L'_0 + L'_1 + \cdots + p + \cdots "1L' + "0L' .  $$

So the corollary follows from Lemma~\ref{lemma:Qomega-countable}.
\end{proof}

\begin{corollary}
Assume that our basic wqo $Q$ is countable. 
W.r.t.\ continuous bi-embeddability 
there are exactly $\omega_1$ many equivalence classes of \qcclos.
\end{corollary}
\begin{proof}
It is easy to see (using the countable ordinals) that the number 
of equivalence classes is at least $\aleph_1$. 

% ***
% ***
% ***
% ***
% ***
% ***

On the other hand,  Corollary~\ref{cor:S-countable} implies that 
 $|\sC_{\alpha }|\le \aleph_0$  for all $\alpha <\omega_1$, so 
 $|\sC_{\omega_1}|\le \aleph_1$. 
\end{proof}

\section{Gödel logics}\label{sec:goedel}

% ***
% ***

Gödel logics are one of the oldest and most interesting families of
many-valued logics. Propositional finite-valued Gödel logics were
introduced by Gödel in \cite{goedel33} to show that intuitionistic logic
does not have a characteristic finite matrix. They
provide the first examples of intermediate logics
(intermediate, that is, in strength between classical and
intuitionistic logics).  Dummett \cite{dummett} was the first to
study infinite valued Gödel logics, axiomatizing the set of
tautologies over infinite truth-value sets by intuitionistic logic
extended by the linearity axiom $(A \limp B) \lor (B \limp A)$. Hence,
infinite-valued propositional Gödel logic is also called
Gödel-Dummett logic or Dummett's \textsc{LC}. In terms of Kripke semantics,
the characteristic linearity axiom picks out those accessibility
relations which are linear orders.  

Quantified propositional Gödel logics and first-order Gödel logics are
natural extensions of the propositional logics introduced by Gödel and
Dummett. For both propositional quantified and first-order Gödel
logics it turns out to be inevitable to consider more complex truth
value sets than the standard unit interval.

Gödel logics occur in a number of different areas of logic
and computer science.  For instance, Dunn and Meyer \cite{DM} pointed
out their relation to relevance logics; Visser \cite{visser} employed
\textsc{LC} in investigations of the provability logic of Heyting
arithmetic; three-valued Gödel logic has been used to model
strong equivalence between logic programs. Furthermore, these logics
have recently received increasing attention, both in terms of
foundational investigations and in terms of applications, as they have
been recognized as one of the most important formalizations of fuzzy
logic \cite{hajek}. 

Perhaps the most surprising fact is that whereas there is only one
infinite-valued propositional Gödel logic, there are infinitely many
different logics at the first-order level
\cite{BaazLeitZach96TCS,Baaz96GOEDEL,Prei02LPAR}.   
In the light of the general result of Scarpellini \cite{scarpellini}
on non-axiomatizability, it is interesting that
some of the infinite-valued Gödel logics belong to the limited class
of recursively enumerable linearly ordered first-order logics
\cite{horn,TT}.

% ***
% ***
% ***
% ***
% ***
% ***
% ***
% ***
% ***
% ***
% ***
% ***
% ***
% ***
% ***
% ***
% ***
% ***
% ***
% ***
% ***
% ***
% ***

Recently a full characterization of axiomatizability of Gödel logics
was given \cite{Prei03PHD}, where also the compactness of the
entailment relation is discussed. But one of the most basic questions
has been left open until now: {\em How many Gödel logics are there?}
Lower bounds to this question have been given in
\cite{Baaz96GOEDEL,Prei02LPAR}, and special subclasses of logics
determined by ordinals have been discussed \cite{ono},
but it was a long open question whether there are only countably many
or uncountably many different Gödel logics. 

% ***

% ***
% ***
% ***
% ***

% ***
% ***
% ***
% ***

% ***

\subsection{Syntax and Semantic}

In the following we fix a relational language  $\LL$ 
of predicate logic with finitely or countably many predicate symbols. 
In addition to the two quantifiers $\forall $ and $\exists$ we 
 use the connectives $\vee$, $\wedge$, $\to$
and the constant $\bot$ (for `false');  other connectives are
introduced as abbreviations, in particular we let
$\lnot \vhi := (\vhi \to \bot)$. 

Originally, Gödel logics have been defined only based on the fixed
truth value set $[0,1]$. But we can fix a (nearly) arbitrary subset
of $[0,1]$ and consider the Gödel logic induced by this truth value
set.

\begin{definition}[Gödel set]
  A Gödel set is any closed set of real numbers, $V \subseteq [0,1]$
  which contains $0$ and $1$.  
\end{definition}

The (propositional) operations on Gödel sets which are used in
defining the semantics of Gödel logics have the property that they are
projecting, i.e.\ that the operation uses one of the arguments (or $1$)
as result: 
\begin{definition} For $a,b\in [0,1]$ let $a\wedge b := \min(a,b)$, 
  $a\vee b:=\max(a,b)$, 
  \[ a \to b := \left\{
    \begin{array}{ll}
      1 & \mbox{if $a\le b$ }\\
      b & \mbox{otherwise}\\
    \end{array}\right.
  \]
The last operation is called `Gödel's implication'.  
Note that 
% ***
$$(a\to b)\  =\ \sup \{\, x: \, (x\wedge a)\, \le\, b\, \};$$
in order theory this is
expressed as `the maps $x\mapsto (a\wedge x)$ and $y\mapsto (a\to y)$
are residuated'.

We define $\lnot a:= (a \to 0)$, so $\lnot 0 = 1$, and $\lnot a = 0$ for all 
$a>0$. 
\end{definition}

% ***
% ***
% ***
% ***
% ***
% ***
% ***
% ***
% ***
% ***
% ***
% ***

The semantics of Gödel logics, with respect to a fixed Gödel set as
truth value set and a fixed relational language~$\LL$ of predicate logic,
is defined using the extended language $\LL^M$, where $M$ is a universe of
objects. 
$\LL^M$ is $\LL$ extended with symbols for
every element of $M$ as constants, so called $M$-symbols. 
These symbols are denoted with the same letters.

\begin{definition}[Semantics of Gödel logic]\label{def:gsemantik}
  Fix a Gödel set $V$ (and a relational language $\LL$). 
  A valuation $v$ into $V$ consists of 
  \begin{enumerate}
  \item  a nonempty set $M = M^v$, the
    `universe' of $v$, 
  \item  for each $k$-ary predicate symbol $P$, 
    a function $P^v : M^k \to V$. 
  \end{enumerate}

  Given a valuation $v$, we can naturally define a value $v(A)$ for
  any closed formula $A$ of $\LL^M$
For atomic formulas $\vhi = P(m_1,\ldots, m_n)$, 
we define $v(\vhi) = P^v( m_1, \ldots, m_n)$, and for
  composite formulas $\vhi$ we define $v(\vhi)$ naturally by: 
  % ***
  \begin{align}
    \val(\bot) &= 0\\
    \val(\vhi\land \psi) &= \min(\val(\vhi),\val(\psi))\\
    \val(\vhi\lor \psi) &= \max(\val(\vhi),\val(\psi))\\
    \val(\vhi\limp \psi) &= \val(\vhi) \to \val(\psi)\\
    \val(\qa x\vhi(x)) &= \inf \{\val(\vhi(m)) \suchthat m\in M\}\\
    \val(\qe x\vhi(x)) &= \sup \{\val(\vhi(m)) \suchthat m\in M\}
  \end{align}
  (Here we use the fact that % ***
our Gödel sets $V$ are \emph{closed} subsets
  of $[0,1]$, in order to be able to interpret $\forall $ and $\exists$ 
  as $\inf$ and $\sup$ in V.)
  
  For any closed formula $\vhi$ and any Gödel set $V$ we let 
  $$ \|\vhi\|_V :=  \inf \{ v(\vhi): \mbox{$v$  a valuation into $V$}\}$$
\end{definition}

\begin{remark} 
  Note that the recursive computation of $v(\vhi)$ depends 
  only on the values $M^v$, $P^v$ and not directly on the 
  set $V$.   Thus, if $V_1\subseteq V_2$ are both Gödel sets, and $v$
  is a valuation into $V_1$, then $v$ can be seen also as a valuation
  into $V_2$, and the values $v(\vhi)$, computed recursively using 
  (1)--(6), 
  do not depend on whether we view $v$ as a $V_1$-valuation or a 
  $V_2$-valuation.

  If $V_1 \subseteq V_2$, there are more valuations into $V_2$ than
  into $V_1$.
  Hence $\|\vhi\|_{V_1} \ge \|\vhi\|_{V_2} $ for all closed $\vhi$.
% ***

  Similarly, for any  map   $h: V_1 \to V_2$, any 
  valuation $v_1$ into $V_1$ induces a valuation $v_2$ 
  into $V_2$ as follows: 
  \[
    M^ {v_1} = M^ {v_2},\quad P^{v_1}(\vec m) = h(P^{v_2}(\vec m)).
  \]
  If  $h: V_1 \to V_2$ is a \smc embedding
  from $V_1$ into $V_2$ which moreover preserves $0$ and $1$,
  and if $v_2$ is the valuation induced by $v_1$ and $h$, then 
  it is easy to verify by induction on the complexity of the
  closed formula $\vhi$ that $v_2(\vhi)=h(v_1(\vhi))$, and hence
  \[ h(\|\vhi\|_{V_1}) \ge \|\vhi\|_{V_2} \]
  for all closed formulas $\vhi$.
\end{remark}

\begin{definition}[Gödel logics based on $V$]
\label{def:goedellogics}
  For a Gödel set $V$ we define the {\em first order Gödel logic $\gdl V$}
as the set of
  all closed formulas of $\LL$ such that $\|\vhi\|_V = 1$.
\end{definition}

From the above remark it is obvious that if $h$ is as above 
or $V_1\subseteq V_2$, the Gödel logic $\gdl{V_2}$ is a subset of $\gdl{V_1}$.

\begin{definition}[Submodel, elementary submodel]
  Let $v_1$, $v_2$ be valuations.  We write $v_1 \subseteq v_2$ 
  ($v_2$ extends $v_1$) iff
  $M^{v_1} \subseteq M^{v_2}$, and for all $k$, all $k$-ary 
  predicate symbols $P$ in $\LL$, we have
  \[ P^{ v_1 } =  P^{ v_2 } \upharpoonright (M^{v_1})^k \]
  or in other words, if $v_1$ and $v_2$ agree on closed atomic
  formulas.  

  We write $v_1 \prec v_2$ if $v_1 \subseteq v_2$ and $v_1(\vhi) = v_2(\vhi)$ 
  for all $\LL^{M^{v_1}}$-formulas $\vhi$. 
\end{definition}

\begin{fact}[downward Löwenheim-Skolem]\label{fact:loewenheim}
  For any valuation $v$ (with $M^v$ infinite) there is a valuation 
  $v' \prec v$  with a  countable universe~$M^{v'} $. 
\end{fact}

% ***
% ***
% ***
% ***

% ***
% ***
% ***
% ***
% ***
% ***
% ***
% ***

% ***
% ***
% ***
% ***
% ***

\begin{definition}
  The only sub-formula of an atomic formula $P$ in $\LL^M$ is  $P$ itself. The
  sub-formulas of $\vhi\star \psi$ for $\star\in\{\limp,\land,\lor\}$ are  
  the subformulas of  $\vhi$ and of $\psi$, together with $\vhi\star  \psi$
  itself.  The sub-formulas of $\qa x\vhi(x)$ and $\qe x\vhi(x)$ with
  respect to a universe $M$   are all subformulas of all 
   $\vhi(m)$ for $m\in M$, together with  $\qa x\vhi(x)$  (or, 
    $\qe x\vhi(x)$, respectively) itself. 

  The set of valuations of sub-formulas of $\vhi$ under a given
  valuation~$\val$ is denoted with
  \[ \subval(\val,\vhi) = \{ \val(\psi) \suchthat \psi 
  \text{ sub-formula of $\vhi$ w.r.t.\ $M^v$}\}\]
\end{definition}

\begin{lemma}\label{lemma:valuation-cut-off}
  Let $\val$ be a valuation with $\val(\vhi)<b<1$
  and $b$ does not occur in $\subval(\val,\vhi)$.
  Let $\val'$ be the valuation with the same universe as $\val$,
  defined by 
  \[ \val'(\psi) = \begin{cases} \val(\psi) &\text{ if $\val(\psi)<b$}\\
                              1       &\text{ otherwise}
                \end{cases}\]
  for atomic subformulas $\psi$  of $\vhi$ w.r.t.~$M^\val$, 
  and arbitrary for all other atomic formulas.
  Then $\val'$ is a valuation and $\val'(\vhi) = \val(\vhi)$.
\end{lemma}

\begin{proof}
  Let $h_b(a)=a$ if $a<b$ and $=1$ otherwise.
  By induction on the complexity of the formula $\psi$
  we can easily show that $\val'(\psi)=h_b(\val(\psi))$ for all
subformulas $\psi$ of $\vhi$ w.r.t.\ $M^\val$.
% ***
\end{proof}

% ***
% ***
% ***
% ***
% ***
% ***
% ***
% ***
% ***
% ***
% ***

% ***
% ***
% ***
% ***
% ***
% ***
% ***

\begin{lemma}\label{lemma:countable-into-cantor}
  Assume that $M\subset \bbR$ is a countable set and $P$ a perfect
  set. Then there is a \smc embedding from $M$ into $P$.
% ***
% ***
% ***
\end{lemma}

In \cite{Prei03PHD} there is a proof of this lemma which was used to
extend the proof of recursive axiomatizability of `standard' Gödel
logics (those with $V=[0,1]$) to Gödel logics with a truth value set
containing a perfect set in the general case. Here we give a 
simple proof.

% ***

\begin{proof}
Since there are uncountable many disjoint sets of the form 
$\bbQ - x :=\{q-x: q\in \bbQ\}$,
there is some $x$ such that $M\cap (\bbQ -x) = \emptyset$, so 
also $(M+x)\cap \bbQ=\emptyset$. 
So we may assume 
% ***
that $M\cap \bbQ=\emptyset$.    
We may also assume $M\subseteq [0,1]$. 

Since $P$ is perfect, we can find an \smc embedding $c$
from the Cantor set $C\subseteq [0,1]$ into $P$.   

Let $i$ be the natural bijection from $2^\omega$ (the set of infinite
$\{0,1\}$-sequences, ordered lexicographically) onto $C$.   $i$ 
is an order preserving homeomorphism. 

For every $m\in M$ let $w(m)\in 2^\omega$ be the binary 
representation of $m$. Since $M$ does not contain any 
dyadic rational numbers, this representation is unique;  moreover,
the map $w$ is smc.   Now $c\circ i \circ  w$ is an smc embedding from $M$
into $P$. 
$$  M  \pfeil{w}  2^\omega \pfeil{i} C  \pfeil{c} P $$

% ***
% ***
% ***
% ***
% ***
% ***

% ***
% ***
% ***
% ***
% ***
% ***
% ***
% ***
% ***
% ***
% ***
% ***
% ***
% ***
% ***
% ***
% ***
% ***
% ***
% ***
% ***
% ***
% ***
% ***
% ***

% ***

% ***
% ***
% ***
% ***
% ***
% ***

% ***
% ***
% ***
% ***
% ***
% ***

% ***
% ***
% ***
% ***
% ***
% ***

% ***
% ***
% ***
% ***

% ***
% ***
% ***
% ***
% ***
% ***
% ***
% ***
% ***
% ***
\end{proof}

\begin{lemma}\label{lemma:gs-id}
  Let $V$ be a truth value set with non-empty perfect kernel $P$, and
  let $W= V \cup [\inf P,1]$, then the logics induced by $V$ and $W$
  are the same, i.e.\ $\gdl{V} = \gdl{W}$.
\end{lemma}

\begin{proof}
% ***
% ***
% ***
% ***
% ***
% ***
As $V\subseteq W$ we have $\gdl{W}\subseteq \gdl V$.
(Cf.\ Remark before Definition~\ref{def:goedellogics}.)

  Now assume that $\val_{W}(\vhi)<1$. Due to Fact~\ref{fact:loewenheim},
  there is a $\val'_W$ such that $M^{v'}$ is countable
  and $\val_W'(\vhi) = \val_W(\vhi)$. The set $M := \subval(v'_W,\vhi)$ has 
  cardinality at most $\aleph_0$, thus there exists a 
  $b\in W$ such that $b\notin M$, $\val'_{W}(\vhi)<b<1$. 
  According to Lemma~\ref{lemma:countable-into-cantor} there is a
  \smc embedding $h$ from $[\inf P, b]\cap (M\cup\{b\})$ into $P$.
  Define  $\val_V(\psi)$ for all atomic subformulas of $\vhi$ as follows: 
  \[ 
    \val_V(\psi) =
      \begin{cases} 
        \val'_W(\psi)     & \text{ if $0<\val_W'(\psi)<\inf P$}\\
        h(\val'_{W}(\psi))& \text{ if $\inf P \le \val_W'(\psi)\le b$}\\
        1                 & \text{ otherwise}
      \end{cases}
  \]
  and $1$ for all other atomic formulas. According to
  Lemma~\ref{lemma:valuation-cut-off} we obtain that 
  \[ \val_V(\vhi) = 
     \begin{cases}
       \val'_W(\psi) < b < 1 & \text{  if $0<\val_W'(\psi)<\inf P$}\\
       h(\val_W'(\vhi)) < h(b) \le 1 &
                    \text{ if $\inf P \le \val_W'(\psi)\le b$}
     \end{cases}
  \]
  thus $\val_V(\vhi) < 1$ and $\gdl V\subseteq \gdl{W}$.
\end{proof}

% ***
% ***
% ***
% ***
% ***
% ***
% ***
% ***
% ***
% ***
% ***
% ***
% ***
% ***
% ***
% ***
% ***
% ***
% ***
% ***
% ***
% ***
% ***
% ***
% ***
% ***
% ***
% ***
% ***
% ***
% ***
% ***
% ***
% ***
% ***
% ***
% ***
% ***
% ***
% ***
% ***

\begin{lemma}\label{lemma:goedelsmc}
  Let $V_1$ and $V_2$ be Gödel sets and $Q=\{0,1\}$ with $0<_Q1$.
  Let $A_1$ and $A_2$ be $Q$-labeled cclos defined by $\dom(A_i)=V_i$, 
  $A_i(0)=A_i(1)=1$ and
  $A_i(x)=0$ otherwise. If $A_1$ is ($Q$-smc-)embeddable into $A_2$,
  then the Gödel logic determined by $V_1$ is a superset of the Gödel
  logic determined by $V_2$.
\end{lemma}

\begin{proof}
  In this case of a very simple labeling the property that $A_1$ is
  embeddable into $A_2$ reduces to the existence of a smc-embedding of
  $V_1$ into $V_2$ preserving $0$ and $1$. According to the Remark
  following Definition~\ref{def:gsemantik} this induces the reverse
  inclusion of the respective Gödel logics.
\end{proof}

\begin{corollary}\label{cor:countable}
  The set of Gödel logics
  \begin{enumerate}
  \item[(a)] is countable 
  \item[(b)] is a (lightface) $\Sigma^1_2$ set
  \item[(c)] is a subset of Gödel's constructible universe $L$.
  \end{enumerate}
\end{corollary}
\begin{proof}
% ***
% ***
% ***
% ***
% ***
% ***
% ***
% ***
% ***
% ***
% ***
% ***
% ***
% ***
% ***
% ***
% ***
(a) First note that the set of countable Gödel logics (i.e.\ those
with countable truth value set), ordered by $\supseteq$, 
is a wqo. 
To see this, assume that $\seqof{\gdl n}{n\in\omega}$ is a sequence of 
countable Gödel logics.
Take the sequence of countable Gödel sets $\seqof{V_n}{n\in\omega}$ 
generating these logics and define the
respective $Q$-labeled cclo (also denoted with $V_n$) with
$Q=\{0,1\}$, $0<_Q1$ and $V_n(0) = V_n(1) = 1$, and $V_n(x) = 0$ otherwise. 
According to Corollary~\ref{cor:sc-wqo} this sequence of $Q$-labeled cclos 
must be good, hence there are numbers $n<m$ such that $V_n$ is 
\smc embeddable into $V_m$.
Then Lemma~\ref{lemma:goedelsmc} implies that 
$\gdl n$ must be a superset of $\gdl m$.
This shows that the original sequence of Gödel logics
$\seqof{\gdl n}{n\in\omega}$ must be good, too.

As each countable Gödel logic is a subset of a fixed countable set (the set 
of all formulas), the family of countable Gödel logics cannot contain a 
copy of $\omega_1$. So by Lemma~\ref{omega1.wqo}, the family of
countable Gödel logics must be countable.

According to Lemma~\ref{lemma:gs-id} any uncountable Gödel logic, i.e.\ 
Gödel logic determined by an uncountable Gödel set, 
such that $0$ is not included in the prefect kernel $P$ of the Gödel set
is completely determined by the countable part
$V\cap[0,\inf P]$. 
So the total number of Gödel logics is at most two times the
number of countable Gödel logics plus $1$ for the logic based on the
full interval, i.e. countable.

% ***
% ***
% ***
% ***
% ***
% ***
% ***

(b) First, note that the set 
$$ \{ (v, \vhi, v(\vhi)):  M^v = \bbN \}$$
is a Borel set, since we can show by induction on the quantifier 
complexity of $\vhi$ that the sets $\{(v, q) : M^v=\bbN, v(\varphi)\ge q\}$
are Borel sets (even of finite rank). 

Next, as set $G$ of formulas is a Gödel logic iff 
\begin{quotation}
 There {\em exists} a closed set $V \subseteq [0,1]$ (say, coded 
as the complement of a sequence of finite intervals) such that: 
\begin{itemize}
\item For every $\vhi\in G$, {\em for every $v$} with $M^v=\bbN$, $v(\vhi)=1$, 
and
\item For every $\vhi\notin G$, there exists $v$ with $M^v=\bbN$, $v(\vhi)<1$.
\end{itemize}
\end{quotation}

(We can restrict our attention to valuations $v$ with $v^M=\bbN$ because of
 Fact~\ref{fact:loewenheim}.)

Counting quantifiers we see that this is a $\Sigma^1_2$ property.

(c) follows from (a) and (b) by the 
Mansfield-Solovay theorem  (see \cite{Mansfield:1970},
\cite[8G.1 and 8G.2]{Moschovakis:1980}). 
\end{proof}

\section*{Questions and future work}

Define % ***
$\omG $ as the smallest ordinal
$\alpha$ such that: For every well-ordered Gödel set $V$ there is a 
well-ordered Gödel set $V'$ of order type $< \alpha$ with $\gdl V=\gdl{V'}$. 

% ***
% ***
% ***
% ***

Define $ \omGCB $ as the smallest ordinal
$\alpha$ such that: For every  Gödel set $V$ there is a 
 Gödel set $V'$ whose Cantor-Bendixson rank is $< \alpha$ with $\gdl V=\gdl{V'}$. 

By Corollary~\ref{cor:countable}, both these ordinals are countable.  
Furthermore,  $\omG  \le \omGCB$. 
% ***
% ***
% ***
% ***
% ***
% ***
% ***
% ***
% ***
% ***
% ***
% ***
% ***
% ***
% ***
% ***
% ***
% ***
% ***
% ***
% ***
% ***
% ***
% ***
It would be interesting to describe the ordinals $\omG$ and $\omGCB$ 
by giving lower and upper estimates in terms of well-known 
closure ordinals, e.g.\ for inductive definitions and related 
reflection principles of set theory.
Are they equal? 
Note that 
$\omega_1^{CK}\le \omG$.

\let\H\hungarian

% \bibliography{bgp}

\end{document}